\documentclass[a4paper,11pt]{article}
\usepackage[utf8]{inputenc}
\usepackage{amsfonts}
\usepackage[english]{babel}
\usepackage{xcolor}
\usepackage[T1]{fontenc}
\usepackage{amsmath, amssymb,epic,graphicx,mathrsfs,enumerate}
\usepackage[all]{xy}
\usepackage{color}
\usepackage{comment}
\usepackage{enumitem}
\usepackage{esint}
\usepackage{hyperref}
\usepackage{amsthm}
\usepackage{latexsym}
\usepackage{epsfig}
\usepackage{soul}
\usepackage{geometry}
\newtheorem{theorem}{Theorem}[section]
\newtheorem{corollary}[theorem]{Corollary}
\newtheorem{lemma}[theorem]{Lemma}
\newtheorem{prop}[theorem]{Proposition}

\title{Existence and compactness of conformal metrics on the plane with unbounded and sign-changing Gaussian curvature}
\author{Chiara Bernardini\\
\small{Università di Padova}\\
\small\texttt{chiara.bernardini@math.unipd.it}}
\date{}

\begin{document}
\maketitle

\begin{abstract}
\noindent We show that the prescribed Gaussian curvature equation in $\mathbb{R}^2$ 
$$-\Delta u= (1-|x|^p) e^{2u},$$
has solutions with prescribed total curvature equal to $\Lambda:=\int_{\mathbb{R}^2}(1-|x|^p)e^{2u}dx\in \mathbb{R}$, if and only if $$p\in(0,2) \qquad \text{and} \qquad (2+p)\pi\le\Lambda<4\pi$$
and prove that such solutions remain compact as $\Lambda\to\bar{\Lambda}\in[(2+p)\pi,4\pi)$, while they produce a spherical blow-up as $\Lambda\uparrow4\pi$.\\



\end{abstract}


\section{Introduction}

In this paper, we study existence and compactness of entire solutions of the equation
\begin{equation}\label{1}
-\Delta u=K e^{2u},\quad\text{in}\,\,\mathbb{R}^2
\end{equation}
where $K\in L^\infty_\text{loc}(\mathbb{R}^2)$ is a given function. Equation \eqref{1} is the prescribed Gaussian curvature equation on $\mathbb{R}^2$. This means that if $u$ satisfies \eqref{1}, then the metric $e^{2u}|dx|^2$ has Gaussian curvature equal to $K$. Equations of this kind also appear in physics, see for example Bebernes and Ederly \cite{BE}, Chanillo and Kiessling \cite{CK} and Kiessling \cite{K}.

The study of equation \eqref{1} started with the work of Liouville \cite{Lio} who considered the case $K=1$; later around the thirties Ahl'fors \cite{A} studied solutions to equation \eqref{1} when $K$ is a negative constant, proving that in this case \eqref{1} does not have any solution in the entire space $\mathbb{R}^2$ (see also the works of Wittich \cite{W}, Sattinger \cite{S} and Oleinik \cite{O}).

The first existence result for solutions to equation \eqref{1} on the entire plane, when the function $K$ is nonpositive and satisfies further growth conditions at infinity, was given by Ni \cite{N} and it was then refined by McOwen \cite{McO} using a weighted Sobolev space approach. Later, a complete classification of all possible solutions to \eqref{1} in some important cases when $K$ is nonpositive was obtained in \cite{cl,CN1,CN2,cl001}.

If $K$ is a positive constant Chen and Li \cite{CL} studied the following problem
\begin{equation}\label{due}
\begin{cases}
-\Delta u=Ke^{2u}\quad\text{on}\,\,\mathbb{R}^2\\ \int_{\mathbb{R}^2}e^{2u}dx<\infty,
\end{cases}
\end{equation}
proving that every solution to \eqref{due} is radially symmetric with respect to some point in $\mathbb{R}^2$. In particular, all such solutions have the form 
\begin{equation}\label{st}
u(x)=\log\frac{2\lambda}{1+\lambda^2|x-x_0|^2}-\log\sqrt{K},
\end{equation}
where $\lambda>0$ and $x_0\in\mathbb{R}^2$. 

If $K(x)$ is positive in some region, under suitable assumptions on the behaviour of $K$ at infinity, existence of solutions to equation \eqref{1} has been intensely studied e.g. by McOwen \cite{McO2}, Aviles \cite{Av} and Cheng and Lin \cite{cl97,cl99,cl001}. 

The study of compactness of solutions to \eqref{1} started with the seminal paper of Brézis and Merle \cite{BM}, which led to a broader study, both in dimension $2$ (see for instance the work of Li and Shafrir \cite{LS}), greater than $2$ (using powers of the Laplacian, or GJMS-operators see e.g. \cite{DM,DR,MS,Mconc,Nd}), or in dimension $1$ (using the $1/2$-Laplacian see e.g. \cite{dalmar,dalmarriv}).

In this paper we will focus on a specific problem that arose recently from a work of Borer, Galimberti and Struwe \cite{BGS} in the context of prescribing Gaussian curvature on $2$-dimensional surfaces. More in detail, they studied a sign-changing prescribed Gaussian curvature problem on a closed, connected Riemann surface $(M,g_0)$ of genus greater than $1$, under the assumption that the prescribed curvature $f_0\le0$ has only non-degenerate maxima $\xi_0$ with $f_0(\xi_0)=0$. In particular, defining $f_\lambda:=f_0+\lambda$ for $\lambda\in\mathbb{R}$, through a mountain-pass technique they investigated the blow-up behavior of ``large''  solutions $u^\lambda$ as $\lambda\downarrow0$. Upon rescaling  they obtain solutions to \eqref{1} with $K$ either constant or  $K(x)=1+(Ax,x)$ where $A=\frac{1}{2}Hess_{f_0}(\xi^{(i)}_\infty)$ (see Theorem 1.4 in \cite{BGS}). More recently, Struwe \cite{S20} obtained a more precise characterisation of this ``bubbling'' behavior: considering a closed surface of genus zero, he proved that all ``bubbles'' are spherical and in fact there are no solutions to equation
\begin{equation}\label{eqS}
-\Delta u=(1+(Ax,x))e^{2u},\quad\text{in}\,\,\mathbb{R}^2
\end{equation}
(where $A$ is a negative definite and symmetric $2\times 2$ matrix) with $u\le C$, such that the induced metric $e^{2u}|dx|^2$ has finite volume, $e^{2u}\in L^1$ and $\int_{\mathbb{R}^2}(1+(Ax,x))e^{2u}dx\in \mathbb{R}$.
In this way he proved that all blow-ups must be spherical also in the higher genus case.

If we do not assume that the prescribed curvature $f_0$ in the problem of Borer, Galimberti and Struwe \cite{BGS} is smooth, but we only require that $f_0\in C^{0,\alpha}$ for some $\alpha\in(0,1]$, then $Hess_{f_0}$ is no more defined and upon rescaling, instead of \eqref{eqS}, one might expect to find solutions to
\begin{equation}\label{7}
-\Delta u= (1-|x|^p) e^{2u}\quad\text{in}\quad\mathbb{R}^2, \qquad \Lambda:=\int_{\mathbb{R}^2}(1-|x|^p)e^{2u}dx<\infty,
\end{equation}
for some $p\in(0,2)$. More precisely, for $p>0$ we define
$$\Lambda_\text{sph}:=4\pi\qquad\text{and}\qquad\Lambda_{*,p}:=(2+p)\pi.$$
(the constant $\Lambda_\text{sph}$ is the total curvature of the sphere $S^2$). It follows from standard arguments the following non-existence result:

\begin{theorem}\label{teo1.1}
Let $p>0$ be fixed. For any $\Lambda\in(-\infty,\Lambda_{*,p})\cup[\Lambda_ \mathrm{sph},+\infty)$ problem \eqref{7} admits no solutions. In particular, for $p\ge2$, problem \eqref{7} admits no solutions.
\end{theorem}

Following an approach of Hyder and Martinazzi \cite{HM20}, and in sharp contrast with the non-existence result of Struwe \cite{S20}, we will show that problem \eqref{7} has solutions for every $p\in(0,2)$ and for suitable $\Lambda$.

\begin{theorem}\label{teo1.2}
Let $p\in(0,2)$ be fixed. Then for every $\Lambda\in(\Lambda_{*,p},\Lambda_ \mathrm{sph})$ there exists a (radially symmetric) solution to problem \eqref{7}. Such solutions have the following asymptotic behavior
\begin{equation}\label{8}
u(x)=-\frac{\Lambda}{2\pi}\log|x|+C+O(|x|^{-\alpha}),\quad\text{as}\,\,|x|\rightarrow\infty,   
\end{equation}
for every $\alpha\in[0,1]$ such that $\alpha<\frac{\Lambda-\Lambda_{*,p}}{\pi}$, and \begin{equation}\label{9}
|\nabla u(x)|=O\left(\frac{1}{|x|}\right), \quad\text{as}\,\, |x|\rightarrow\infty.
\end{equation}
\end{theorem}

Observe that Theorem \ref{teo1.1} and Theorem \ref{teo1.2} do not cover the case $\Lambda=\Lambda_{*,p}$. In this case (see Proposition \ref{2.3} and Lemma \ref{2.5}) relation \eqref{8} degenerates to 
$$-\frac{\Lambda_{*,p}+o(1)}{2\pi}\log|x|\le u(x)\le-\frac{ \Lambda _{*,p}}{2\pi} \log|x|+O(1), \quad \text{as}\,\,|x|\to+\infty,$$ 
which is compatible with the integrability of $(1-|x|^p)e^{2u}$. We will study the case $\Lambda=\Lambda_{*,p}$ from the point of view of compactness, namely we will show that solutions to \eqref{7} are compact for $\Lambda$ away from $\Lambda_{\mathrm{sph}}$, and blow up spherically at the origin as $\Lambda\uparrow \Lambda_{\mathrm{sph}}$.

\begin{theorem}\label{teo1.3}
Fix $p\in(0,2)$, let $(u_k)$ be a sequence of solutions to \eqref{7} with $\Lambda=\Lambda_k\in[\Lambda_{*,p},\Lambda_\mathrm{sph})$ and $\Lambda_k\rightarrow\bar{\Lambda}\in[\Lambda_{*,p},\Lambda_\mathrm{sph})$. Then, up to a subsequence, $u_k\rightarrow \bar{u}$ locally uniformly, where $\bar{u}$ is a solution to \eqref{7} with $\Lambda=\bar{\Lambda}$.\\
\indent Moreover, choosing $\Lambda_k\downarrow\Lambda_{*,p}$ and $u_k$ given by Theorem \ref{teo1.2}, we obtain that \eqref{7} has a solution $u$ also for $\Lambda=\Lambda_{*,p}$ and we have
\begin{equation}\label{10}
u(x)\le-\frac{\Lambda_{*,p}}{2\pi}\log|x|-\left(1+o(1)\right)\log\log|x|,\quad \text{as}\,\, |x|\rightarrow+\infty
\end{equation}
and 
\begin{equation}\label{10a}
|\nabla u(x)|=O\left(\frac{1}{|x|}\right),\quad \text{as}\,\, |x|\rightarrow+ \infty.\end{equation}
\end{theorem}

The proof of Theorem \ref{teo1.3} relies on uniform controls of the integral of the Gaussian curvature at infinity. This is particularly subtle as $\Lambda_k\downarrow \Lambda_{*,p}$ since in this case  $Ke^{2u_k}$ is a priori no better than uniformly $L^1$ at infinity and we could have loss of negative curvature at infinity. This is ruled out with an argument based on the Kelvin transform (see Lemma 5.4). 

Theorem \ref{teo1.3} strongly uses the lack of scale invariance of equation \eqref{7}. Indeed for the constant curvature case \eqref{due}, solutions given in \eqref{st} can blow up ``spherically''. We will show that this is also the case when $\Lambda\uparrow\Lambda_\text{sph}$.

\begin{theorem}\label{teo1.4}
Fix $p\in(0,2)$, let $(u_k)$ be a sequence of radial solutions to \eqref{7} with $\Lambda=\Lambda_k\uparrow\Lambda_\mathrm{sph}$ as $k\to+\infty$. Then 
$$(1-|x|^p)e^{2u_k}\rightharpoonup\Lambda_\mathrm{sph}\delta_0,\quad\mathrm{as}\,\,k\to+\infty,$$
weakly in the sense of measures. Moreover, setting
$$\mu_k:=2e^{-u_k(0)}$$
and
$$\eta_k(x):=u_k(\mu_k x)-u_k(0)+\log2,$$
we have $$\eta_k(x)\xrightarrow[k\rightarrow\infty]{}\log\frac{2}{1+|x|^2} \quad\text{in}\,\,\,C^1_{\text{loc}}(\mathbb{R}^2).$$
\end{theorem}

Notice that all solutions to problem \eqref{7} are radially symmetric about the origin and monotone decreasing, this follows by Corollary 1 in the work of Naito \cite{Nai} (to prove Corollary 1 Naito uses an approach  based on the maximum principle in unbounded domains together with the method of moving planes).\\

The outline of the paper is the following. In Section 2 we prove some regularity results for solution to \eqref{7}. In Section 3, thanks to some estimates on the asymptotic behavior of solutions at infinity and to a Pohozaev-type identity, we prove Theorem \ref{teo1.1}. In Section 4, we use a variational approach due to Chang and Chen \cite{CC} to prove Proposition \ref{3.1}, which, together with a blow-up argument will allow us to prove the existence part of Theorem \ref{teo1.2}. 
Section 5 is devoted to the proof of Theorem \ref{teo1.3}, which is based on a blow-up analysis, the Kelvin transform, and on quantization of the total curvature (using Theorem 2 in \cite{M}). Finally, the proof of Theorem \ref{teo1.4} in Section 6 is based on the fact that if $\Lambda\uparrow\Lambda_\text{sph}$, from Theorem \ref{teo1.1}, we could have only two cases: loss of curvature at infinity, or loss of compactness; from Lemma \ref{5.1} we get that the second case occurs.\\

After the completion of the work the author learned that a result of Cheng and Lin (see Theorem 1.1 in \cite{cl99}) implies our Theorem \ref{teo1.2}. Their elegant proof is based on a Moser-Trudinger inequality in weighted Sobolev spaces. Our approach to Theorem \ref{teo1.2}, only based on ODE methods, is more elementary and will be the basis for the compactness Theorem \ref{teo1.3}, therefore we left the statement of Theorem \ref{teo1.2} for completeness.


\section{Regularity of solutions}

Let $p>0$ be fixed. First of all, we prove some regularity results for solutions to equation
\begin{equation}\label{eqreg}
-\Delta u=(1-|x|^p)e^{2u}\quad\text{in}\,\,\mathbb{R}^2,
\end{equation}
assuming that $u\in L^1_{loc}(\mathbb{R}^2)$ and $(1-|x|^p)e^{2u}\in L^1(\mathbb{R}^2)$.\\

In the following, if $\Omega\subseteq\mathbb{R}^2$ is an open set and $s\in\mathbb{R}$ $s\ge0$, we will denote
$$C^s(\Omega):=\left\{u\in C^{\lfloor s\rfloor}(\Omega)\,\,\big|\,\,D^{\lfloor s\rfloor}u\in C^{0,s-\lfloor s \rfloor}(\Omega)\right\},$$
and 
$$C^s_{loc}(\mathbb{R}^2):=\left\{u\in C^0(\mathbb{R}^2)\,\,\big|\,\,u|_{\Omega}\in C^s(\Omega)\,\,\text{for every}\,\,\Omega\subset\subset\mathbb{R}^2\right\}.$$

\begin{prop}\label{W}
Let $u$ be a solution to \eqref{eqreg}. Then $u\in W^{2,r}_{loc} (\mathbb{R}^2)$ for $1<r<\infty$.
\end{prop}

\proof
Following the proof of Theorem 2.1 in \cite{HMM}, first we prove that $e^{2u}\in L^q_{loc}(\mathbb{R}^2)$ for any $q\ge1$. Let $q\ge1$ be fixed, take $\varepsilon=\varepsilon(q)$ such that $q<\frac{\pi}{2\varepsilon}$, we can find two functions $f_1$ and $f_2$ such that $(1-|x|^p)e^{2u}=f_1+f_2$ and 
$$f_1\in L^1(\mathbb{R}^2)\cap L^\infty(\mathbb{R}^2), \qquad \|f_2\|_{\mathrm{L}^1(\mathbb{R}^2)}<\varepsilon.$$
We define the following 
$$u_i(x):=\frac{1}{2\pi}\int_{\mathbb{R}^2}\log\left(\frac{|y|}{|x-y|}\right)f_i(y)dy,\quad i=1,2,$$
and $$u_3:=u-u_1-u_2.$$ In this way, $u_3$ is harmonic and hence $u_3\in C^\infty(\mathbb{R}^2)$. Differentiating $u_1$ we obtain 
$$\nabla u_1(x)=-\frac{1}{2\pi}\int_{\mathbb{R}^2} \frac{x-y}{|x-y|^2}f_1(y)dy$$
using the relation $|x||y|\Bigl|\frac{x}{|x|^2}-\frac{y}{|y|^2}\Bigr|=|x-y|$ it's easy to prove that $\nabla u_1$ is continuous, hence $u_1\in C^{1}(\mathbb{R}^2)$. Concerning $u_2$ we have
\begin{align*}
\int_{B_R}e^{8qu_2}dx&=\int_{B_R}\exp{\left(\int_{\mathbb{R}^2}\frac{8q\|f_2\|}{2\pi}\log\left(\frac{|y|}{|x-y|}\right)\frac{f_2(y)}{\|f_2\|}dy\right)}dx\\
&\le\int_{B_R}\int_{\mathbb{R}^2}\exp\left(\frac{8q\|f_2\|}{2\pi}\log\left(\frac{|y|}{|x-y|}\right)\right)\frac{f_2(y)}{\|f_2\|}dy\,dx\\
&=\frac{1}{\|f_2\|}\int_{\mathbb{R}^2}f_2(y)\int_{B_R}\left(\frac{|y|}{|x-y|}\right)^{\frac{4q\|f_2\|}{\pi}}dx\,dy\le C,
\end{align*}
using Holder's inequality, we can conclude that $e^{2u}\in L^q_{loc} (\mathbb{R}^2)$ for any $q\ge1$.

By assumption $p>0$, so $(1-|x|^p)\in L^r_{loc}(\mathbb{R}^2)$ for every $1\le r\le\infty$; it follows that $-\Delta u=(1-|x|^p)e^{2u}\in L^r_{loc}(\mathbb{R}^2)$ for each $1\le r\le\infty$. By elliptic estimates (refer to \cite[Theorem 9.11]{GT}), we have $$u\in W^{2,r}_{loc}(\mathbb{R}^2),\quad\text{for every}\,\,r\in(1,\infty),$$
and by the Morrey-Sobolev embedding we get
$$u\in C^{1,\alpha}_{loc}(\mathbb{R}^2),\quad\text{for}\,\,\alpha\in(0,1].$$
\endproof

\begin{prop}\label{reg}
Let $p>0$ be fixed and $u$ be a solution of equation \eqref{eqreg}. If $p\not\in\mathbb{N}$ then $u\in C^\infty(\mathbb{R}^2\setminus\{0\})\cap  C^{p+2}_{loc}(\mathbb{R}^2)$, if $p-1\in2\mathbb{N}$ then $u\in C^\infty(\mathbb{R} ^2\setminus\{0\})\cap C^{s+2}_{loc}(\mathbb{R}^2)$ for $s<p$ and if $p\in2\mathbb{N}$ then $u\in C^\infty(\mathbb{R}^2)$.
\end{prop}

\proof
From Proposition \ref{W} we have that $u\in W^{2,r}_{loc}(\mathbb{R}^2)\hookrightarrow C^{1,\alpha}_{loc}(\mathbb{R}^2)$ for $1<r<\infty$ and $\alpha\in(0,1]$. Since $1-|x|^p\in C^\infty(\mathbb{R}^2 \setminus\{0\})$ if $p\not\in2\mathbb{N}$, and belongs to $C^\infty(\mathbb{R}^2)$ if $p\in2\mathbb{N}$, by bootstrapping regularity we get that $u\in C^\infty(\mathbb{R}^2\setminus\{0\})$ if $p\not\in2\mathbb{N}$ and $u\in C^\infty(\mathbb{R}^2)$ if $p\in2\mathbb{N}$. Moreover, we can verify that for $p>0$
$$1-|x|^p\in\begin{cases}
C^p(\mathbb{R}^2)\quad\text{if}\,\,p\not\in\mathbb{N}\\
C^s_{loc}(\mathbb{R}^2)\quad\text{for}\,\,s<p,\,\,\text{if}\,\,p-1\in2\mathbb{N}\\
C^\infty(\mathbb{R}^2) \quad\text{if}\,\,p\in2\mathbb{N}
\end{cases}.$$
Using Schauder estimates and bootstrapping regularity, we get that if $p\not\in\mathbb{N}$ then $u\in C^{p+2}_{loc}(\mathbb{R}^2)$ and if $p-1\in2 \mathbb{N}$ then $u\in C^{s+2}_{loc}(\mathbb{R}^2)$ for $s<p$.
Hence we can conclude that if $p\not\in\mathbb{N}$ then $u\in C^\infty(\mathbb{R}^2\setminus\{0\})\cap  C^{p+2}_{loc}(\mathbb{R}^2)$, if $p-1\in2\mathbb{N}$ then $u\in C^\infty(\mathbb{R} ^2\setminus\{0\})\cap C^{s+2}_{loc}(\mathbb{R}^2)$ for $s<p$ and if $p\in2\mathbb{N}$ then $u\in C^\infty(\mathbb{R}^2)$.
\endproof

All solutions to \eqref{eqreg} are in fact normal solutions, namely if $u$ is a solution to \eqref{eqreg} then $u$ solves the integral equation
\begin{equation}\label{2}
u(x)=\frac{1}{2\pi}\int_{\mathbb{R}^2}\log\left(\frac{|y|}{|x-y|}\right)(1-|y|^p)e^{2u(y)}dy+c,
\end{equation}
where $c\in\mathbb{R}$ (for a detailed proof of this fact see the proof of Theorem 2.1 in \cite{cl97}). Moreover, if we have more integrability (namely $\log(|\cdot|)(1-|y|^p)e^{2u}\in L^1(\mathbb{R}^2)$), equation \eqref{2} is equivalent to
\begin{equation}\label{3}
u(x)=\frac{1}{2\pi}\int_{\mathbb{R}^2}\log\left(\frac{1}{|x-y|}\right)(1-|y|^p)e^{2u(y)}dy+c',   
\end{equation}
where $c'\in\mathbb{R}$.


\section{Non-existence result}
First of all, we can observe that the case $\Lambda<0$ is not possible (this follows similar to the proof of Theorem 1 in \cite{Mar}).

\begin{prop}\label{lambda>0}
Let $p>0$ be fixed, if $u$ is a solution to \eqref{7} then $\Lambda\ge0$.
\end{prop}

\proof 
Assume by contradiction that $\Lambda<0$, then there exists $r_0>0$ such that
$$\int_{B_r}\Delta u \, dx\ge0,\quad\forall r\ge r_0,$$
hence
$$\int_{\partial B_r}\frac{\partial u}{\partial\nu}d\sigma(x)\ge0,\quad\forall r\ge r_0.$$
It follows that $\fint\limits_{\partial B_r}u\,d\sigma$ is an increasing function for $r\ge r_0$, and consequently also $\exp\left(\fint_{\partial B_r}u\, d\sigma\right)$ is increasing for $r\ge r_0$. By Jensen inequality we get
$$\exp\left(2\fint_{\partial B_r}u\,d\sigma\right)\le\fint_{\partial B_r}e^{2u}d\sigma.$$
It follows that $\fint_{\partial B_r}e^{2u}d\sigma$ must be increasing for $r\ge r_0$, hence integrating $\int_{\mathbb{R}^2}e^{2u}=+\infty$. This leads to a contradiction since in this way $\Lambda$ cannot be finite.
\endproof

In order to prove Theorem \ref{teo1.1} we need the following Pohozaev-type identity and some asymptotic estimates at infinity. 

\begin{prop}\label{2.1}
In the case $K(x)=1-|x|^p$, let $u$ be a solution to the integral equation
\begin{equation}\label{17}
u(x)=\frac{1}{2\pi}\int_{\mathbb{R}^2}\log\left(\frac{|y|}{|x-y|}\right)K(y)e^{2u(y)}dy+c,
\end{equation}
for some $c\in\mathbb{R}$, with $Ke^{2u}\in L^1(\mathbb{R}^2)$ and $|\cdot|^pe^{2u}\in L^1(\mathbb{R}^2)$. If 
\begin{equation}\label{18}
\lim\limits_{R\rightarrow\infty}R^{2+p}\max\limits_{|x|=R}e^{2u(x)}=0
\end{equation}
then denoting $\Lambda=\int_{\mathbb{R}^2}K(x)e^{2u(x)}dx$, we have
\begin{equation}\label{19}
\frac{\Lambda}{\Lambda_\mathrm{sph}}(\Lambda-\Lambda_\mathrm{sph})=-\frac{p}{2}\int\limits_{\mathbb{R}^2}|x|^p e^{2u(x)}dx.
\end{equation}
\end{prop}

\proof
In the spirit of the proof of Theorem 2.1 in \cite{XU},  differentiating equation \eqref{17} (by Proposition \ref{reg} $u$ is sufficiently regular) and multiplying by $x$, we obtain
$$\langle x,\nabla u(x)\rangle=-\frac{1}{2\pi}\int_{\mathbb{R}^2}\frac{\langle x,x-y\rangle}{|x-y|^2}K(y)e^{2u(y)}dy.$$
Multiplying both sides of the previous one by $K(x)e^{2u(x)}$ and integrating over the ball $B_R(0)$ for $R>0$, we have
\begin{align*}
\int_{B_R}K(x)&e^{2u(x)}\langle x,\nabla u(x)\rangle dx=\\
&=\int_{B_R}K(x)e^{2u(x)}\left[-\frac{1}{2\pi}\int_{\mathbb{R}^2}\frac{\langle x,x-y\rangle}{|x-y|^2}K(y)e^{2u(y)}dy\right]dx.  
\end{align*}
Integrating by parts the left-hand side:
\begin{align*}
\int_{B_R}&K(x)e^{2u(x)}\langle x,\nabla u(x)\rangle dx
=\frac{1}{2}\int_{B_R}K(x)\langle x,\nabla e^{2u(x)}\rangle dx=\\
&=-\int_{B_R}\left(K(x)+\frac{1}{2}\langle x,\nabla K(x)\rangle\right)e^{2u(x) }dx +\frac{R}{2}\int_{\partial B_R}K(x)e^{2u(x)}d\sigma= \\
&=-\int_{B_R}\left(K(x)-\frac{p}{2}|x|^p\right)e^{2u(x)}dx+\frac{R}{2}\int_{\partial B_R}(1-|x|^p)e^{2u(x)}d\sigma.
\end{align*}
It's easy to see that
$$-\int_{B_R}\left(K(x)-\frac{p}{2}|x|^p\right)e^{2u(x)}dx\xrightarrow[R\rightarrow\infty]{}-\Lambda+\frac{p}{2}\int_{\mathbb{R}^2}|x|^p e^{2u(x)}dx$$
and concerning the boundary term, we have
\begin{align*}
\frac{R}{2}\int_{\partial B_R}&(1-|x|^p)e^{2u(x)}d\sigma\le \frac{R}{2}\max\limits_{|x|=R}e^{2u(x)}\int_{\partial B_R}(1-|x|^p)d\sigma=\\
&=\pi R^2(1-R^p)\max\limits_{|x|=R}e^{2u(x)}=
\pi R^{p+2}\left(\frac{1}{R^p}-1\right)\max\limits_{|x|=R}e^{2u(x)}
\end{align*}
using \eqref{18} it goes to 0 if $R\rightarrow+\infty$.
Regarding the right-hand side, we have
\begin{align*}
\int_{B_R}K(x)&e^{2u(x)}\left[-\frac{1}{2\pi}\int_{\mathbb{R}^2}\frac{\langle x,x-y\rangle}{|x-y|^2}K(y)e^{2u(y)}dy\right]dx=\\
&=\frac{1}{2}\int_{B_R}K(x)e^{2u(x)}\left[-\frac{1}{2\pi}\int_{\mathbb{R}^2}K(y)e^{2u(y)}dy\right]dx\\
&\quad+\frac{1}{2}\int_{B_R}K(x)e^{2u(x)}\left[-\frac{1}{2\pi}\int_{\mathbb{R}^2}\frac{\langle x+y,x-y\rangle}{|x-y|^2}K(y)e^{2u(y)}dy\right]dx.
\end{align*}
It follows immediately that 
$$\frac{1}{2}\int_{B_R}K(x)e^{2u(x)}\left[-\frac{1}{2\pi}\int_{\mathbb{R}^2}K(y)e^{2u(y)}dy\right]dx\xrightarrow[R\rightarrow+\infty]{}-\frac{1}{4\pi}\Lambda^2$$
while, by dominated convergence, the second term goes to $$-\frac{1}{4\pi}\int_{\mathbb{R}^2}\int_{\mathbb{R}^2}\frac{\langle x+y,x-y\rangle}{|x-y|^2}K(x)e^{2u(x)}K(y)e^{2u(y)}dy\,dx$$
which, changing variables, it's equal to 0. Finally combining we obtain \eqref{19}.
\endproof

In order to study the asymptotic behavior of solution to \eqref{1}, a useful trick is the Kelvin transform. We have the following result.

\begin{prop}\label{2.2}
Let $u$ be a solution to \eqref{1} with $K(x)=1-|x|^p$ and $Ke^{2u}\in L^1(\mathbb{R}^2)$. Then the Kelvin transform of $u$, namely the function
\begin{equation}\label{20} \tilde{u}(x)=u\left(\frac{x}{|x|^2}\right)-\alpha\log{|x|},\qquad \text{for}\quad x\not=0,
\end{equation}
where $\alpha:=\frac{1}{2\pi}\Lambda$, satisfies
$$\tilde{u}(x)=\frac{1}{2\pi}\int_{\mathbb{R}^2}\log\left(\frac{1}{|x-y|}\right)K\left(\frac{y}{|y|^2}\right)\frac{e^{2\tilde{u}(y)}}{|y|^{4-2\alpha}}dy+c,$$
and hence $\tilde{u}$ is a solution to $$-\Delta\tilde{u}(x)=K\left(\frac{x}{|x|^2}\right)\frac{e^{2\tilde{u}}}{|x|^{4-2\alpha}}.$$
\end{prop}

\proof
Following the proof of Proposition 2.2 in \cite{HM20}, using \eqref{2} and then changing variables, we have
\begin{align*}
\tilde{u}(x)&=\frac{1}{2\pi}\int_{\mathbb{R}^2}\log\left(\frac{|y|}{\bigl|\frac{x}{|x|^2}-y\bigr|}\right)K(y)e^{2u(y)}dy-\frac{1}{2\pi}\log |x|\int_{\mathbb{R}^2}K(y)e^{2u(y)}dy+c=\\
&=\frac{1}{2\pi}\int_{\mathbb{R}^2}\log\left(\frac{|y|}{|x|\Bigl|\frac{x}{|x|^2}-y\Bigr|}\right)K(y)e^{2u(y)}dy+c=\\
&=\frac{1}{2\pi}\int_{\mathbb{R}^2}\log\left(\frac{1}{|x||y|\Bigl|\frac{x}{|x|^2}-\frac{y}{|y|^2}\Bigr|}\right)K\left(\frac{y}{|y|^2}\right)\frac{e^{2\tilde{u}(y)}}{|y|^{4-2\alpha}}dy+c=\\
&=\frac{1}{2\pi}\int_{\mathbb{R}^2}\log\left(\frac{1}{|x-y|}\right)K\left(\frac{y}{|y|^2}\right)\frac{e^{2\tilde{u}(y)}}{|y|^{4-2\alpha}}dy+c
\end{align*}
where in the last equality we have used that $|x||y|\Bigl|\frac{x}{|x|^2}-\frac{y}{|y|^2}\Bigr|=|x-y|$.
\endproof

\subsection{Asymptotic behavior of solutions}

\begin{prop}\label{2.3}
Let $u$ be a solution to the integral equation
$$u(x)=\frac{1}{2\pi}\int_{\mathbb{R}^2}\log\left(\frac{|y|}{|x-y|}\right)K(y)e^{2u(y)}dy+c,$$
where $c\in\mathbb{R}$, $K(y)\le0$ for $|y|\ge R_0>0$ and $Ke^{2u}\in L^1$. Then we have
\begin{equation}\label{23}
u(x)\le-\frac{\Lambda}{2\pi}\log|x|+O(1),\qquad\text{as}\quad|x| \rightarrow\infty. \end{equation}
\end{prop}

\proof
Following the proof of Lemma 2.3 in \cite{HM20} and adapting it to the two-dimensional case, we choose $x$ such that $|x|\ge 2R_0$ (assuming $R_0\ge 2$) and consider $\mathbb{R}^2=A_1\cup A_2\cup A_3$ where
$$A_1=B_{\frac{|x|}{2}}(x),\qquad A_2=B_{R_0}(0), \qquad A_3=\mathbb{R}^2\setminus(A_1\cup A_2).$$
It's easy to see that $$\int_{A_1}\log\left(\frac{|y|}{|x-y|}\right)K(y)e^{2u(y)}dy\le 0$$
because if $y\in A_1$ we have $K(y)\le0$ and $\log\left(\frac{|y|}{|x-y|}\right)\ge0$. 
For $y\in A_2$ we have $\log\left(\frac{|y|}{|x-y|}\right)=-\log|x|+O(1)$ as $|x|\rightarrow+\infty$, hence
$$\int_{A_2}\log\left(\frac{|y|}{|x-y|}\right)K(y)e^{2u(y)}dy=-\log|x|\int_{A_2}Ke^{2u}dy+O(1),\quad\text{as}\,\,|x|\to+\infty.$$
If $y\in A_3$ we have $|x-y|\le|x|+|y|\le|x||y|$, hence $K(y)\log\left(\frac{|y|}{|x-y|}\right)\le K(y)\log\left(\frac{1}{|x|}\right)$, so we obtain
$$\int_{A_3}\log\left(\frac{|y|}{|x-y|}\right)K(y)e^{2u(y)}dy\le -\log|x|\int_{A_3}Ke^{2u}dy.$$
Finally we have $$u(x)\le\frac{1}{2\pi}\left(-\log|x|\int_{A_2\cup A_3}Ke^{2u}dy\right)+O(1)\le-\frac{\Lambda}{2\pi}\log|x|+O(1),$$
using the fact that $\int_{A_2\cup A_3}Ke^{2u}\ge\int_{\mathbb{R}^2} Ke^{2u}=\Lambda$ since $K\le0$ in $A_1$.
\endproof

\begin{corollary}\label{2.4}
If $p\in(0,2)$, there exist no solutions to \eqref{7} for $\Lambda\ge4\pi$.
\end{corollary}

\proof
Assume that $u$ solves \eqref{7} for some $\Lambda\ge4\pi$. By Proposition \ref{2.3} $u$ satisfies \eqref{23}, so hypothesis \eqref{18} in Proposition \ref{2.1} is verified and therefore from \eqref{19} we must have $\Lambda<4\pi$, a contradiction.
\endproof

\begin{lemma}\label{2.5}
Fix $p>0$ and let $u$ be a solution to \eqref{7}. Then we have 
$$\Lambda\ge\Lambda_{*,p}$$
and 
\begin{equation}\label{24} u(x)=-\frac{\Lambda+o(1)}{2\pi}\log|x|,\quad\mathrm{as}\,\,|x|\rightarrow+\infty \end{equation}
\end{lemma}

\proof
Following the proof of Lemma 2.5 in \cite{HM20}, we first prove the asymptotic estimate \eqref{24}. Writing $u=u_1+u_2$, where $$u_2(x)=-\frac{1}{2\pi}\int_{B_1(x)}\log\left(\frac{1}{|x-y|}\right)|y|^p e^{2u(y)}dy,$$
we obtain
\begin{equation}\label{stimau1}
u_1(x)=-\frac{\Lambda+o(1)}{2\pi}\log|x|,\quad\mathrm{as}\,\,|x|\rightarrow+\infty, \end{equation}
(see the Appendix for a detailed proof of the estimate \eqref{stimau1}).
Let's consider $R\gg1$ and $|x|\ge R+1$, define 
$$h(R):=\frac{1}{2\pi}\int_{B_R^c}|y|^p e^{2u}dy,$$
it's easy to see that $h(R)\to 0$ as $R\to+\infty$, so we can write $h(R)=o_R(1)$. We have
$$-u_2(x)=\int_{B_R^c}h(R)\log\left(\frac{1}{|x-y|}\right)\chi_{|x-y|\le1}d\mu(y),\quad d\mu(y)=\frac{|y|^p e^{2u}}{\int_{B_R^c}|y|^p e^{2u}dy}dy.$$
By Jensen's inequality and Fubini's theorem we get 
$$\int_{R+1<|x|<2R}e^{-2u_2}dx\le\int_{B_R^c}\int_{R+1<|x|<2R}\left(1+\frac{1}{|x-y|^{2h(R)}}\right)dx\,d\mu(y)\le CR^2.$$
Hence using Holder's inequality and the previous one
\begin{equation}\label{form}
R^2\approx\int_{R+1<|x|<2R}e^{u_2}e^{-u_2}dx\le CR\left(\int_{R+1<|x|<2R} e^{2u_2}dx\right)^{1/2}.
\end{equation}
If we assume by contradiction that $\frac{\Lambda}{\pi}\le p$, then  $|y|^p e^{2u_1}\ge\frac{1}{|y|}$ for $|y|$ sufficiently large, therefore
$$o_R(1)=\int_{R+1<|x|<2R}|x|^p e^{2u_1}e^{2u_2}dx\gtrsim
\frac{1}{R}\int_{R+1<|x|<2R}e^{2u_2}dx$$ which contradicts \eqref{form}.
Therefore we must have $\frac{\Lambda}{\pi}>p$, from this it follows that $|y|^p e^{2u_1}\le C$ on $\mathbb{R}^2$, using the fact that $u_2\le0$ we have $|y|^p e^{2u_1}e^{2u_2}<C$, and then
$$|u_2(x)|\le C\int_{B_1(x)}\log\frac{1}{|x-y|}dy\le C$$
this prove \eqref{24}. Finally, if $\Lambda<\Lambda_{*,p}$ then $(1-|x|^p)e^{2u}\not\in L^1(\mathbb{R}^2)$, hence it must be $\Lambda\ge\Lambda_{*,p}$.
\endproof

\begin{lemma}\label{2.7}
Let $u$ be a solution to \eqref{7} with $\Lambda=\Lambda_{*,p}=(2+p)\pi$ and $\tilde{u}$ its Kelvin transform (as in \eqref{20}). Then
\begin{equation}\label{26} \lim\limits_{x\to 0}\tilde{u}(x)=-\infty.\end{equation}
and
\begin{equation}\label{27} \lim\limits_{x\to 0}\Delta\tilde{u}(x) =+\infty.\end{equation}
\end{lemma}

\proof
In this case, we have
\begin{equation}\label{25}
\tilde{u}(x)=u\left(\frac{x}{|x|^2}\right)-\left(1+\frac{p}{2}\right)\log{|x|},\quad \text{for}\,\, x\not=0,
\end{equation}
and $\tilde{u}$ satisfies
\begin{equation}\label{28}
\tilde{u}(x)=\frac{1}{2\pi}\int_{\mathbb{R}^2}\log\left(\frac{1}{|x-y|}\right)\left(1-\frac{1}{|y|^p}\right)\frac{e^{2\tilde{u}(y)}}{|y|^{2-p}}dy+c,\quad \text{for}\,\,x\not=0.
\end{equation}
We follow the proof of Lemma 2.7 in \cite{HM20}. Using Proposition \ref{2.3}, which gives us a bound from above, we have $\sup_{B_1}\tilde{u}<+\infty$. Since $\tilde{u}\le C$ in $B_1$ and $u$ is continuous, using \eqref{25}, we get
\begin{equation}\label{29}
\left(1-\frac{1}{|y|^p}\right)\frac{e^{2\tilde{u}(y)}}{|y|^{2-p}}
\le\frac{C}{|y|^4}\qquad \text{on}\quad B_1^c.
\end{equation}
Considering again equation \eqref{28}, we have
\begin{align*}
\tilde{u}(x)=&-\frac{1}{2\pi}\int\limits_{B_1}\log\left(\frac{1}{|x-y|}\right)\frac{e^{2\tilde{u}(y)}}{|y|^2}dy+\frac{1}{2\pi}\int\limits_{B_1}\log\left(\frac{1}{|x-y|}\right)\frac{e^{2\tilde{u}(y)}}{|y|^{2-p}}dy+\\
&+\frac{1}{2\pi}\int\limits_{B_1^c}\log\left(\frac{1}{|x-y|}\right)\left(1-\frac{1}{|y|^p}\right)\frac{e^{2\tilde{u}(y)}}{|y|^{2-p}}dy+c.
\end{align*}
If $0<|x|<1$ the second and the third term of above are $O(1)$, therefore
\begin{equation}\label{30}
\tilde{u}(x)=-\frac{1}{2\pi}\int\limits_{B_1}\log\left(\frac{1}{|x-y|}\right)\frac{e^{2\tilde{u}(y)}}{|y|^2}dy+O(1),\qquad \text{if}\,\,\,0<|x|<1.
\end{equation}
Assume by contradiction that $\tilde{u}(x_k)=O(1)$  for a sequence $x_k\to 0$, applying Lemma 2.6 of \cite{HM20} to \eqref{30} we have
$$\int_{B_1}\log\left(\frac{1}{|y|}\right)\frac{e^{2\tilde{u}(y)}}{|y|^2}dy<+\infty$$
and changing variables
\begin{equation}\label{30a} 
\int_{B_1^c}\log(|y|)|y|^p e^{2u(y)} dy<+\infty.
\end{equation}
Since $\Lambda=\Lambda_{*,p}$, as $|x|\to\infty$ we have
\begin{align*}
\int_{|y|\le\sqrt{|x|}}&(1-|y|^p)e^{2u(y)}dy=\Lambda_{*,p}-\int_{|y|>\sqrt{|x|}}(1-|y|^p)e^{2u(y)}dy=\\
&=\Lambda_{*,p}+O\left(\frac{1}{\log|x|}\int_{|y|>\sqrt{|x|}}\log(|y|)(1-|y|^p)e^{2u(y)}dy\right)=\\
&=\Lambda_{*,p}+O\left(\frac{1}{\log|x|}\right) \qquad\text{as}\quad |x|\to+\infty \end{align*}
where in the last equality we used \eqref{30a}. Using equation \eqref{3} for $|x|\gg1$
$$u(x)=\frac{1}{2\pi}\left(\int_{|y|\le\sqrt{|x|}}+\int_{\sqrt{|x|}\le|y|\le2|x|}+\int_{|y|\ge 2|x|}\right)\log\left(\frac{1}{|x-y|}\right)(1-|y|^p)e^{2u(y)}dy+c.$$
From the previous estimate, as $|x|\to\infty$, 
\begin{align*}
\int_{|y|\le\sqrt{|x|}}&\log\left(\frac{1}{|x-y|}\right)(1-|y|^p)e^{2u(y)}dy=\\
&=(-\log|x|+O(1))\left(\Lambda_{*,p}+O\left(\frac{1}{\log|x|}\right)\right)=-\Lambda_{*,p}\,\log|x|+O(1).
\end{align*}
Concerning the second integral, as $|x|\to\infty$ we get
\begin{align*}
&\int_{\sqrt{|x|}\le|y|\le2|x|}\log\left(\frac{1}{|x-y|}\right)(1-|y|^p)e^{2u(y)}dy\ge\int_{B_1(x)}\log\left(\frac{1}{|x-y|}\right)(1-|y|^p)e^{2u(y)}dy\\
&\quad=-\int_{B_1(x)}\log\left(\frac{1}{|x-y|}\right)\big|1-|y|^p\big|e^{2u}dy\ge-\frac{1}{|x|^2}\int_{B_1(x)}\log\left(\frac{1}{|x-y|}\right)dy=O\left(\frac{1}{|x|^2}\right), \end{align*}
and for the third integral, using \eqref{30a}, we have
$$\frac{1}{2\pi}\int_{|y|\ge2|x|}\left(\frac{1}{|x-y|}\right)(1-|y|^p)e^{2u(y)}dy= O(1).$$
Therefore, for $|x|\to \infty$ 
$$u(x)\ge-\frac{\Lambda_{*,p}}{2\pi}\log|x|+O(1)$$
but this means that $|\cdot|^p e^{2u}\not\in L^1(\mathbb{R}^2)$, which is a contradiction. Hence \eqref{26} is proven.\\
\indent In this case, $\tilde{u}$ is a solution to 
$$-\Delta\tilde{u}(x)=\left(1-\frac{1}{|x|^p}\right) \frac{e^{2\tilde{u}(x)}} {|x|^{2-p}},$$
using \eqref{25} and \eqref{24} we obtain that $\lim\limits_{x\to0}(-\Delta\tilde{u}(x))=-\infty$, which proves \eqref{27}.
\endproof

\begin{prop}\label{2.8}
If $p\ge2$, there exists no solution to \eqref{7}.
\end{prop}

\proof
Assume by contradiction that for some $p\ge 2$ there exists a solution $u$ of \eqref{7}, then by Lemma \ref{2.5} we must have $\Lambda\ge\Lambda_{*,p}$ and moreover $\Lambda_{*,p}\ge4\pi$ since $p\ge2$.
If $\Lambda>\Lambda_{*,p}$, using Lemma \ref{2.3} we have $u(x)\le-\frac{\Lambda}{2\pi}\log|x|+c$ for $|x|\to+\infty$. In this way $u$ satisfies \eqref{18} and hence by Proposition \ref{2.1} follows that
$$\Lambda_\text{sph}\le\Lambda_{*,p}<\Lambda<\Lambda_\text{sph},$$
which is a contradiction.
If $\Lambda=\Lambda_{*,p}$, using \eqref{25} and \eqref{26} we observe that \eqref{18} is satisfied, and we proceed as in the previous case.
\endproof

\proof[Proof of Theorem \ref{teo1.1}] 
Theorem \ref{teo1.1} is proven, we have just to combine Corollary \ref{2.4}, Lemma \ref{2.5} and Proposition \ref{2.8}.
\endproof

If $\Lambda=\Lambda_{*,p}$ we obtain a sharper version of \eqref{26}.

\begin{lemma}\label{2.9}
Fix $p\in(0,2)$, let $u$ be a solution to \eqref{7} with $\Lambda=\Lambda_{*,p}$. Then we have $$\limsup\limits_{|x|\to+\infty} \frac{u(x)+\left(1+\frac{p}{2}\right)\log|x|}{\log\log|x|}=-1.$$
\end{lemma}

\proof
We refer to the proof of Lemma 2.9 in \cite{HM20}. Let $\tilde{u}$ be defined as in \eqref{25}, $u$ is a radial solution (this follows from Corollary 1 in \cite{Nai}) 
therefore also $\tilde{u}$ is radially symmetric. By Lemma \ref{2.7} we have $$\lim\limits_{r\to0}\tilde{u}(r)=-\infty,\qquad\lim\limits_{r\to0}\Delta\tilde{u}(r)=+\infty,$$
so there exists $\delta>0$ such that $\tilde{u}$ is monotone increasing in $B_\delta(0)$. Using this and \eqref{30}, we estimate for $|x|\to0$ and get
\begin{align*}
-\tilde{u}(x)&\ge\frac{1}{2\pi}\int_{2|x|\le|y|<1}\log\left(\frac{1}{|x-y|}\right)\frac{e^{2\tilde{u}(y)}}{|y|^2}dy+O(1)\\
&=\frac{1}{2\pi}\int_{2|x|\le|y|\le\delta}\log\left(\frac{1}{|y|}\right)\frac{e^{2\tilde{u}(y)}}{|y|^2}dy+\frac{1}{2\pi}\int_{\delta\le|y|<1}\log\left(\frac{1}{|y|}\right)\frac{e^{2\tilde{u}(y)}}{|y|^2}dy+O(1)\\
&\ge\frac{e^{2\tilde{u}(x)}}{2\pi}\int_{2|x|\le|y|\le\delta}\log\left(\frac{1}{|y|}\right)\frac{dy}{|y|^2}+O(1)\\
&=e^{2\tilde{u}(x)}\int_{2|x|}^\delta\frac{\log\frac{1}{\rho}}{\rho}d\rho+O(1), \quad\text{as}\,\, |x|\to0.
\end{align*}
Hence we have 
$$-\tilde{u}(x)+O(1)\ge \frac{e^{2\tilde{u}(x)}}{2}\left(\log\left(\frac{1}{2|x|}\right) \right)^2, \quad\text{as}\,\, |x|\to0.$$
Now taking the logarithm and rearranging
$$\limsup\limits_{x\to0}\frac{\tilde{u}(x)}{\log\log\left(\frac{1}{|x|}\right)}\le-1.$$
We prove that this $\limsup$ is equal to $-1$. Assume by contradiction, that the previous $\limsup$ is less than $-1$; hence it must exist $\varepsilon>0$ such that  $$\tilde{u}(x)\le-\left(1+\frac{\varepsilon}{2}\right)\log\log \frac{1}{|x|}$$
for $|x|$ small. Therefore, recalling \eqref{30}, for $|x|$ small we have
$$-\tilde{u}(x)\le C\int\limits_{B_1}\log\left(\frac{1}{|x-y|}\right)\frac{dy}{|y|^2|\log|y||^{2+\varepsilon}}+O(1).$$
We can split $\int\limits_{B_1}\log\left(\frac{1}{|x-y|}\right)\frac{dy}{|y|^2|\log|y||^{2+\varepsilon}}$ into $I_1+I_2+I_3$ where
$$I_i=\int\limits_{A_i}\log\left(\frac{1}{|x-y|}\right)\frac{dy}{|y|^2|\log|y||^{2+\varepsilon}}$$
where
$$A_1=B_{\frac{|x|}{2}},\quad A_2=B_{2|x|}\setminus B_{\frac{|x|}{2}}\,\, \text{and}\,\, A_3=B_1\setminus B_{2|x|}.$$
Concerning $I_1$, we observe that if $y\in B_{\frac{|x|}{2}}$ we have $\log\left(\frac{1}{|x-y|}\right)\sim\log\left(\frac{1}{|x|}\right)$ as $|x|\to0$ and $\int_{B_{|x|/2}}\frac{1}{|y|^2|\log|y||^{2+\varepsilon}}dy=\big|\log\frac{|x|}{2}\big|^{-1-\varepsilon}$, hence 
$$I_1\le\frac{C}{|\log|x||^\varepsilon}.$$
Regarding $I_2$ we have 
$$I_2\le\frac{C}{|x|^2\big|\log|x|\big|^{2+\varepsilon}}\int_{B_{2|x|}}\log\left(\frac{1}{|y|}\right)dy=\frac{C}{\big|\log|x|\big|^{1+\varepsilon}}.$$
Finally it's easy to prove that $I_3\le C$. 
We have obtained a contradiction to the fact that $-\tilde{u}(x)\to+\infty$ as $|x|\to0$.
\endproof


\section{Existence result}

The proof of the existence part of Theorem \ref{teo1.2} will be based on a blow-up argument as done in \cite{HM20}, but first we need the following result, which we will prove using a variational approach due to Chang and Chen \cite{CC}.

\begin{prop}\label{3.1}
Let $p\in(0,2)$ be fixed. For every $\Lambda\in(0,\Lambda_\mathrm{sph})$ and $\lambda>0$, there exists a radial solution $u_\lambda$ to 
\begin{equation}\label{34}
-\Delta u_\lambda=(\lambda-|x|^p)e^{-|x|^2} e^{2u_\lambda} \quad\mathrm{on}\,\,\mathbb{R}^2,
\end{equation}
such that
\begin{equation}\label{35}
\Lambda=\int_{\mathbb{R}^2}(\lambda-|x|^p)e^{-|x|^2} e^{2u_\lambda}dx.
\end{equation}
Moreover, such solution $u_\lambda$ solves the integral equation
\begin{equation}\label{ulam}
u_\lambda(x)=\frac{1}{2\pi}\int_{\mathbb{R}^2}\log\left(\frac{1}{|x-y|}\right)(\lambda-|y|^p)e^{-|y|^2} e^{2u_\lambda}dy+c_\lambda
\end{equation}
for some constant $c_\lambda\in\mathbb{R}$.
\end{prop}

\proof
Let's denote $K_\lambda(x)=(\lambda-|x|^p)e^{-|x|^2}$, we prove that there exists at least one radial solution to 
\begin{equation}\label{ulambda}
-\Delta u_\lambda=K_\lambda(x) e^{2u_\lambda}\quad\text{on}\,\,\mathbb{R}^2.
\end{equation}
Following the proof of Theorem 2.1 in \cite{CC}, we identify each point in $\mathbb{R}^2$ with a point on $S^2$ through the stereographic projection $\Pi:S^2\rightarrow\mathbb{R}^2$, take $K(x):=K_\lambda(x)$ and $\mu=1-\frac{\Lambda}{\Lambda_\text{sph}}$. 
A solution to \eqref{ulambda} is of the form 
$$u_\lambda=w\circ\Pi^{-1}+(1-\mu)\eta_0,$$
where $\eta_0(x)=\log\left(\frac{2}{1+|x|^2}\right)$, $w=u+c$ where $u$ minimizes a certain functional defined on the set of functions
$$\left\{v\in H^1_\mathrm{rad}(S^2)\,\bigg|\,\int K_\lambda(x)e^{2v}dV>0\right\}$$
and $c$ is a suitable constant such that $\int K_\lambda(x) e^{2w}dV=(1-\mu)\Lambda_\text{sph}$. By construction \eqref{35} holds 
$$\int K_\lambda(x)e^{2u_\lambda}dx=\int K_\lambda(x)e^{2w}dV=(1-\mu)\Lambda _\text{sph}=\Lambda.$$
In order to prove \eqref{ulam}, consider
$$-\Delta_{g_0}w+(1-\mu)=(K_\lambda\circ\Pi)e^{-2\mu(\eta_0\circ\Pi)}e^{2w}.$$
Since $(K_\lambda\circ\Pi)e^{-2\mu(\eta_0\circ\Pi)} \in L^\infty(S^2)$ and $e^{2w}\in L^q(S^2)$ for every $q\in[1,\infty)$, by elliptic estimates we get $w\in C^{1,\alpha}(S^2)$ for $\alpha\in(0,1)$. Therefore $w$ is continuous in $S=(0,0,-1)$ and hence
$$u_\lambda(x)=(1-\mu)\eta_0(x)+w(S)+o(1)=\frac{\Lambda}{2\pi}\log|x|+C+o(1),\quad\text{as}\,\,|x|\to+\infty.$$
Defining $$v_\lambda:=\frac{1}{2\pi}\int_{\mathbb{R}^2}\log\left(\frac{1}{|x-y|}\right)K_\lambda(y)e^{2u_\lambda(y)}dy$$
and $h_\lambda:=u_\lambda-v_\lambda$, we observe that 
$$\Delta h_\lambda=0, \qquad h_\lambda(x)=O(\log|x|)\quad\text{as}\,\,|x|\to+ \infty$$
by Liouville's theorem $h_\lambda$ must be constant.
\endproof

If $w$ is a radial function belonging to $C^2(\mathbb{R}^2)$ and $0\le r<R$, using the divergence theorem, we have the following identity
\begin{equation}\label{37}
w(r)-w(R)=\int_r^R\frac{1}{2\pi t}\int_{B_t}-\Delta w dx dt.
\end{equation}

Let $u_\lambda$ be a radial solution to \eqref{34} given by Proposition \ref{3.1}.

\begin{lemma}\label{3.2}
For every $\lambda>0$ we have $u_\lambda(x)\downarrow-\infty$ as $|x|\to+\infty$.
\end{lemma}

\proof
Let's consider the function  $$r\longmapsto\int_{B_r}-\Delta u_\lambda(x)dx=\int_{B_r}(\lambda-|x|^p)e^{-|x|^2}e^{2u_\lambda(x)}dx$$
we observe that it's increasing on $[0,\lambda^{1/p}]$ and decreasing to $\Lambda>0$ on $[\lambda^{1/p},+\infty)$, so it follows that it's positive for every $r>0$. Hence by \eqref{37} $u_\lambda$ is a decreasing function of $|x|$ and using Proposition \ref{2.3} we conclude that $u_\lambda\to-\infty$ as $|x|\to+\infty$.
\endproof

\begin{lemma}\label{3.3}
$\lambda e^{2u_\lambda(0)}\to +\infty$ as $\lambda\downarrow0$.
\end{lemma}

\proof
We assume that $\lambda e^{2u_\lambda(0)}\le C$ as $\lambda\downarrow0$, we have
\begin{align*}
\Lambda&=\int_{\mathbb{R}^2}(\lambda-|x|^p)e^{-|x|^2} e^{2u_\lambda}dx\le \int_{B_{\lambda^{1/p}}}(\lambda-|x|^p)e^{-|x|^2} e^{2u_\lambda}dx\\
&\le\int_{B_{\lambda^{1/p}}}\lambda e^{-|x|^2} e^{2u_\lambda(0)}dx\xrightarrow[\lambda\to 0]{}0,
\end{align*}
contradiction.
\endproof

Now we define $$\eta_\lambda(x):=u_\lambda(r_\lambda x)-u_\lambda(0)$$
where $r_\lambda$ are non-negative and defined in such a way that 
$$\lambda\, r_\lambda^2\, e^{2u_\lambda(0)}=1,$$
by Lemma \ref{3.3} we have  $r_\lambda\to0$ as $\lambda\downarrow0$, moreover $\eta_\lambda(0)=0$ and $\eta_\lambda\le0$ (since $u_\lambda$ is a radial decreasing function).
A basic calculation shows that $\eta_\lambda$ is a solution to the equation
$$-\Delta\eta_\lambda=\left(1-\frac{r_\lambda^p |x|^p}{\lambda}\right) e^{-r_\lambda^2|x|^2} e^{2\eta_\lambda}$$
such that
\begin{equation}\label{39}
\Lambda=\int_{\mathbb{R}^2}\left(1-\frac{r_\lambda^p |x|^p}{\lambda}\right) e^{-r_\lambda^2|x|^2} e^{2\eta_\lambda}dx.
\end{equation}
We have 
\begin{equation}\label{39a}
0<\Lambda<\int_{B_{\frac{\lambda^{1/p}}{r_\lambda}}}e^{-r_\lambda^2|x|^2}e^{2\eta_\lambda}dx\le \big|B_{\frac{\lambda^{1/p}}{r_\lambda}}\big|
\end{equation}
where in the second inequality we have used that $\left(1-\frac{r_\lambda^p |x|^p}{\lambda}\right)<0$ in $B^c_{\frac{\lambda^{1/p}}{r_\lambda}}$; instead in the last inequality we employ the fact that since $\eta_\lambda\le0$ then $e^{-r_\lambda^2|x|^2} e^{2\eta_\lambda}\le1$.
Finally from \eqref{39a} follows $$\limsup\limits_{\lambda\to0}\frac{r_\lambda^p}{\lambda}<+\infty.$$

Since $\eta_\lambda(0)=0$, using ODE theory, we have that, up to a subsequence,
\begin{equation}\label{43}
\eta_\lambda\xrightarrow[\lambda\to0]{}\eta\qquad\text{in}\,\,C^2_{loc}(\mathbb{R}^2) \end{equation}
where the function $\eta$ satisfies the following equation in $\mathbb{R}^2$
$$-\Delta\eta=(1-\mu|x|^p)e^{2\eta}, \qquad \mu:=\lim\limits_{\lambda\to0}\frac{r_\lambda^p}{\lambda}\in[0,+\infty).$$
At this stage we do not know if $\mu>0$ and if $\int_{\mathbb{R}^2}(1-|x|^p)e^{2\eta}= \Lambda$.

\begin{lemma}\label{3.5}
If $\mu=0$ then $e^{2\eta}\in L^1(\mathbb{R}^2)$.
\end{lemma}

\proof
From Lemma \ref{3.2} and \eqref{43}, we have that $\eta$ is decreasing.
Since $\eta$ satisfies $-\Delta\eta= e^{2\eta}$, then $\Delta\eta\le0$ must be increasing. We have $\lim_{r\to\infty}\Delta\eta(r)=:c_0\in[-\infty,0]$, if $c_0=0$ then $\lim_{r\to\infty}e^{2\eta(r)}=0$ and so $e^{2\eta}\in L^1(\mathbb{R}^2)$; if $c_0<0$ then $\eta(r)\lesssim-r^2$, and hence $e^{2\eta}\in L^1(\mathbb{R}^2)$.
\endproof

\begin{lemma}\label{3.6}
For every $\Lambda\in(\Lambda_{*,p},\Lambda_\mathrm{sph})$ we have $\mu>0$.
\end{lemma}

\proof 
Assume by contradiction that $\mu=0$, then $\eta$ is a solution to 
\begin{equation}\label{11}
-\Delta\eta=e^{2\eta}\qquad\text{in}\quad\mathbb{R}^2
\end{equation}
where $e^{2\eta}\in L^1(\mathbb{R}^2)$ from Lemma \ref{3.5}. Chen and Li in \cite{CL} proved that every solution to \eqref{11} with finite total Gaussian curvature is a standard one, namely assumes the form 
$$\eta(x)=\log\left(\frac{2\lambda}{1+\lambda^2|x-x_0|^2}\right),$$
for some $\lambda>0$ and $x_0\in\mathbb{R}^2$. Hence all solutions are radially symmetric with respect to some point $x_0\in\mathbb{R}^2$. Since $\eta$ is spherical, we have
$$\int_{\mathbb{R}^2} e^{2\eta}dx=|S^2|=\Lambda_\text{sph}.$$
Moreover, by assumption $\Lambda<\Lambda_\text{sph}$, so we can fix $R_0>0$ such that $$\Lambda<\int\limits_{B_{R_0}} e^{2\eta}dx\,\,.$$ 
Recalling that $r_\lambda\rightarrow 0$ and $\frac{r_\lambda^p}{\lambda}\rightarrow0$ as $\lambda\downarrow0$, one can find $\lambda_0$ (which depends on $R_0$) such that
\begin{equation}\label{45}
\int\limits_{B_{R_0}}\left(1-\frac{r_\lambda^p}{\lambda}|x|^p\right)e^{-r_\lambda^2 |x|^2}e^{2\eta_\lambda}dx\ge\Lambda,\quad \text{for}\quad\lambda\in(0,\lambda_0).
\end{equation}
Defining the function 
$$\Gamma_\lambda(t):=\int\limits_{B_t}\left(1-\frac{r_\lambda^p}{\lambda}|x|^p\right)e^{-r_\lambda^2 |x|^2}e^{2\eta_\lambda}dx,$$
we can observe that $\Gamma_\lambda(0)=0$, $\Gamma_\lambda$ is monotone increasing on the interval $\left[0,\frac{\lambda^{1/p}}{r_\lambda}\right]$ and it decreases to $\Lambda$ on $\left[\frac{\lambda^{1/p}}{r_\lambda},+\infty\right)$. 
From \eqref{45} we have 
\begin{equation}\label{46}
\Gamma_\lambda(t)\ge\Lambda, \quad\text{for}\,\, t\ge R_0\,\,\text{and}\,\,\lambda\in(0,\lambda_0).
\end{equation}
Integrating from $R_0$ to $r\ge R_0$ we have
$$-\int_{R_0}^r\frac{\Gamma_\lambda(t)}{2\pi t}dt\le-\int_{R_0}^r\frac{\Lambda}{2\pi t}dt=-\frac{\Lambda}{2\pi}\log{\frac{r}{R_0}}=-\left(1+\frac{p}{2}+\delta\right)\log{\frac{r}{R_0}}$$
where $\delta>0$ is such that $\Lambda_{*,p}+2\delta\pi=\Lambda$, hence by \eqref{37} we have
\begin{align*}
\eta_\lambda(r)&
\le\eta_\lambda(R_0)-\left(1+\frac{p}{2}+\delta\right)\log{\frac{r}{R_0}}=\\
&=C(R_0)-\left(1+\frac{p}{2}+\delta\right)\log{r}, \qquad r\ge R_0.
\end{align*}
This implies that 
\begin{equation}\label{48}
\lim\limits_{R\rightarrow\infty}\lim\limits_{\lambda\rightarrow 0}\int_{B_R^c}(1+|x|^p)e^{2\eta_\lambda}dx=0.
\end{equation} 
We can split
$$\Lambda=\int_{B_R}\left(1-\frac{r_\lambda^p|x|^p}{\lambda}\right) e^{-r^2_\lambda |x|^2} e^{2\eta_\lambda}dx+\int_{B_R^c}\left(1-\frac{r_\lambda^p|x|^p}{\lambda}\right) e^{-r^2_\lambda |x|^2} e^{2\eta_\lambda}dx$$
by uniform convergence the first term goes to $\int_{B_R}e^{2\eta}dx$ as $\lambda\to0$, regarding the second term we have
$$\bigg|\int_{B_R^c}\left(1-\frac{r_\lambda^p|x|^p}{\lambda}\right) e^{-r^2_\lambda |x|^2} e^{2\eta_\lambda}dx\bigg|\le\int_{B_R^c}(1+|x|^p)e^{2\eta_\lambda}dx,$$
so taking the limit as $\lambda\to 0$ and using \eqref{48}, we obtain 
$$\Lambda=\int_{B_R}e^{2\eta}dx+o_R(1),\qquad\text{as}\,\,R\to+\infty$$
hence it follows that 
$$\Lambda=\int_{\mathbb{R}^2}\left(1-\frac{r_\lambda^p|x|^p}{\lambda}\right) e^{-r^2_\lambda |x|^2} e^{2\eta_\lambda}dx\xrightarrow [\lambda\downarrow0]{}\int_{\mathbb{R}^2}e^{2\eta}dx=\Lambda_\text{sph}$$
which is absurd, hence we must have $\mu>0$.
\endproof

\proof[Proof of the existence part of Theorem \ref{teo1.2}] 
Since, by Lemma \ref{3.6},\, $\frac{r_\lambda^p}{\lambda}\xrightarrow[\lambda\to0]{}\mu>0$, choosing $R=\left(\frac{4}{\mu}\right)^{1/p}$, we have, for $\lambda$ sufficiently small
\begin{equation}\label{48.5}
1-\frac{r_\lambda^p}{\lambda}|x|^p\le1-\frac{\mu}{2}|x|^p\le-\frac{\mu}{4}|x|^p \qquad\text{for}\quad |x|\ge R.
\end{equation}
Therefore, by \eqref{48.5} and \eqref{39}, we obtain
\begin{equation}\label{49}
\begin{split}
\int_{B_R^c}\frac{\mu}{4}|x|^p e^{-r_\lambda^2 |x|^2}e^{2\eta_\lambda}dx&\le-\int_{B_R^c}\left(1-\frac{r_\lambda^p}{\lambda}|x|^p\right) e^{-r_\lambda^2|x|^2}e^{2\eta_\lambda}dx=\\
&=\int_{B_R}\left(1-\frac{r_\lambda^p}{\lambda}|x|^p\right) e^{-r_\lambda^2|x|^2}e^{2\eta_\lambda}dx-\Lambda\le C
\end{split}
\end{equation}
where we used that $\eta_\lambda\le0$. From the fact that the integrand in $B_R$ is uniformly bounded, it follows that 
$$\int_{\mathbb{R}^2}(1+|x|^p)e^{-r_\lambda^2|x|^2}e^{2\eta_\lambda}dx\le C.$$
From \eqref{49} we also have
$$\int_{B_R^c}e^{-r_\lambda^2|x|^2}e^{2\eta_\lambda}dx\le\frac{C}{R^p}\xrightarrow[R\rightarrow+\infty]{} 0$$
uniformly in $\lambda$, and hence
\begin{equation}\label{51}
\lim\limits_{\lambda\to0}\int_{\mathbb{R}^2}e^{-r_\lambda^2 |x|^2} e^{2\eta_\lambda}dx=\int_{\mathbb{R}^2}e^{2\eta}dx.
\end{equation}
Using Fatou's lemma
\begin{equation}\label{52}
\int_{\mathbb{R}^2}|x|^p e^{2\eta}dx\le\lim\limits_{\lambda\to0} \int_{\mathbb{R}^2}|x|^p e^{-r_\lambda^2|x|^2} e^{2\eta_\lambda}dx.
\end{equation}
From \eqref{51} and \eqref{52} we obtain
$$\int_{\mathbb{R}^2}(1-\mu|x|^p)e^{2\eta}dx\ge\Lambda.$$
Now we are going to prove that the previous inequality is actually an equality. Since $\frac{r_\lambda^p}{\lambda}\to\mu>0$ then 
$$\frac{\lambda^{\frac{1}{p}}}{r_\lambda}\to\frac{1}{\mu^{\frac{1}{p}}}>0,$$
proceeding as in the proof of Lemma \ref{3.6}, for $R_0=2\mu^{-\frac{1}{p}}$ and $\lambda_0$ sufficiently small, relation \eqref{46} holds and from it \eqref{48} follows. Finally, \eqref{48} implies that
$$\int_{\mathbb{R}^2}(1-\mu|x|^p)e^{2\eta}dx=\Lambda.$$
Defining $$u(x):=\eta(\rho x)+\log\rho,\qquad\rho:=\mu^{-\frac{1}{p}}$$ 
we get the desired solution to \eqref{7}.
\endproof

\subsection{Asymptotic behaviour}

\proof[Proof of \eqref{8}]
Let $u$ be a solution to \eqref{7} given by Theorem \ref{teo1.2} and $\tilde{u}$ its Kelvin transform as defined in \eqref{20}. We have 
$$-\Delta\tilde{u}(x)=\left(1-\frac{1}{|x|^p}\right)\frac{e^{2\tilde{u}(x)}}{|x|^{4-\frac{\Lambda}{\pi}}}=O\left(\frac{1}{|x|^{4+p-\Lambda/\pi}}\right),\qquad \text{as}\,\,|x|\to0.$$
Since for $\Lambda>\Lambda_{*,p}$ 
$$4+p-\frac{\Lambda}{\pi}=2-\frac{\Lambda-\Lambda_{*,p}}{\pi}<2$$
we obtain 
$$-\Delta\tilde{u}\in L^q_{loc}(\mathbb{R}^2) \quad\text{for}\,\,1\le q<\frac{1}{1-\frac{\Lambda-\Lambda_{*,p}}{2\pi}},$$
hence by elliptic estimates we have  $$\tilde{u}\in\mathrm{W}^{2,q}_\mathrm{loc}(\mathbb{R}^2)\quad \text{for}\,\,1\le q<\frac{1}{1-\frac{\Lambda-\Lambda_{*,p}}{2\pi}}.$$
Now, by the Morrey-Sobolev embedding we obtain that 
$$\tilde{u}\in C_{loc}^{0,\alpha}(\mathbb{R}^2)$$ for $\alpha\in[0,1]$ such that $\alpha<\frac{\Lambda-\Lambda_{*,p}}{\pi}$. From this \eqref{8} follows.
\endproof

\proof[Proof of \eqref{9}]
Let $u$ be a solution to \eqref{7}, then $u$ satisfies the integral equation
$$u(x)=\frac{1}{2\pi}\int_{\mathbb{R}^2}\log\left(\frac{|y|}{|x-y|}\right)(1-|y|^p)e^{2u(y)}dy+c,$$
differentiating under the integral sign, we obtain
$$|\nabla u(x)|=
O\left(\int_{\mathbb{R}^2}\frac{1}{|x-y|}(1+|y|^p)e^{2u(y)}dy\right)$$
since $\Lambda>\Lambda_{*,p}$, using \eqref{8} as $|x|\to\infty$ we get
$$(1+|x|^p)e^{2u(x)}\le\frac{C}{1+|x|^{2+\delta}}$$
for $\delta>0$. Hence, for $|x|$ large, we have
\begin{align*}
|\nabla u(x)|&\le C\left(\int_{B_{|x|/2}}+\int_{B_{2|x|}\setminus B_{|x|/2}}+\int_{B_{2|x|}^c}\right)\frac{1}{|x-y|}\,\frac{1}{1+|y|^{2+\delta}}dy\\
&\le \frac{C}{|x|}+\frac{C}{|x|^{2+\delta}}\int_{B_{2|x|}\setminus B_{|x|/2}}\frac{1}{|x-y|}dy\le \frac{C}{|x|}.
\end{align*}
\endproof


\section{Proof of Theorem \ref{teo1.3}}

Let $p\in(0,2)$ be fixed and $(u_k)$ be a sequence of solutions to \eqref{7} with $\Lambda=\Lambda_k \in [\Lambda_{*,p},\Lambda_\text{sph})$, hence every solution $u_k$ solves the integral equation
\begin{equation}\label{58}
u_k(x)=\frac{1}{2\pi}\int_{\mathbb{R}^2}\log\left(\frac{|y|}{|x-y|}\right)(1-|y|^p) e^{2u_k(y)}dy+c_k.
\end{equation}
Assuming that
\begin{equation}\label{59}
\Lambda_k:=\int_{\mathbb{R}^2}(1-|x|^p)e^{2u_k(x)}dx\rightarrow\bar{\Lambda}\in[\Lambda_{*,p},\Lambda_\text{sph}),
\end{equation}
we want to prove that $u_k\rightarrow\bar{u}$ (up to a subsequence) uniformly locally in $\mathbb{R}^2$, where $\bar{u}$ is a radial solution to \eqref{7} with $\Lambda=\bar{\Lambda}$.

With the same procedure used in the proof of Lemma \ref{3.2}, we can prove that $u_k$ is radially decreasing.

\begin{lemma}\label{4.2}
We have $u_k(0)\ge C$ where $C$ depends only on $\inf\limits_k \Lambda_k$.
\end{lemma}

\proof Using the fact that $u_k$ is radially decreasing, we have
\begin{align*}
\Lambda_k&=\int_{\mathbb{R}^2}(1-|x|^p)e^{2u_k(x)}dx\le
\int_{B_1}(1-|x|^p)e^{2u_k(x)}dx\\
&\le\int_{B_1}e^{2u_k(x)}dx\le\pi e^{2u_k(0)},
\end{align*}
therefore 
$$u_k(0)\ge\frac{1}{2}\log{\frac{\Lambda_k}{\pi}}\ge\log{\sqrt{p+2}}>0.$$
\endproof

Since $\Lambda_k\in[\Lambda_{*,p},\Lambda_\text{sph})$, using Proposition \ref{2.1} (which can be applied thanks to Proposition \ref{2.3} if $\Lambda_k\in (\Lambda_{*,p}, \Lambda_\text{sph})$ and thanks to Lemma \ref{2.9} if $\Lambda_k=\Lambda_{*,p}$) we have the following Pohozaev identity:
\begin{equation}\label{60}
\frac{\Lambda_k}{\Lambda_\text{sph}}(\Lambda_k-\Lambda_\text{sph})=-\frac{p}{2}\int_{\mathbb{R}^2}|x|^p e^{2u_k(x)}dx.
\end{equation}
Hence we have
\begin{equation}\label{61}
\int_{\mathbb{R}^2}e^{2u_k}dx=\Lambda_k-\frac{2\Lambda_k}{p\Lambda_\text{sph}}(\Lambda_k-\Lambda_{\text{sph}})
\end{equation}
and taking the limit for $k\rightarrow+\infty$
\begin{equation}\label{62}
\lim_{k\rightarrow+\infty}\int_{\mathbb{R}^2}e^{2u_k}dx=\bar{\Lambda}-\frac{2\bar{\Lambda}}{p\Lambda_\text{sph}}(\bar{\Lambda}-\Lambda_{\text{sph}}).
\end{equation}

\begin{lemma}\label{4.3}
$$\limsup\limits_{k\rightarrow+\infty}{u_k(0)}<\infty$$
\end{lemma}

\proof 
From \eqref{59}, \eqref{60} and \eqref{61} we have
\begin{equation}\label{63}
\limsup_{k\to+\infty}\int_{\mathbb{R}^2}(1+|x|^p)e^{2u_k}dx<+\infty.
\end{equation}
Differentiating \eqref{58}, integrating over $B_1(0)$ and using Fubini's theorem, we get
$$\int_{B_1(0)}|\nabla u_k|dx\le C \int_{\mathbb{R}^2}\left(\int_{B_1(0)}\frac{1} {|x-y|}dx\right)(1+|y|^p)e^{2u_k(y)}dy\le C$$
where in the last inequality we used \eqref{63}. Assume by contradiction that (up to a subsequence) $u_k(0)\rightarrow+\infty$ as $k\rightarrow+\infty$, by Theorem 2 in \cite{M} (the two-dimensional case was first studied by Brezis and Merle in \cite{BM} and by Li and Shafrir in \cite{LS}) we have that $u_k\rightarrow-\infty$ locally uniformly in $B_1(0)\setminus\{0\}$.
Moreover, we have quantization of the total curvature, namely
$$\lim\limits_{k\to+\infty}\int_{B_r(0)}(1-|x|^p)e^{2u_k}dx=4\pi,$$
where $r\in(0,1)$ is fixed. 
The blow-up at the origin is spherical, namely there exists a sequence of positive numbers 
$$\mu_k:=2e^{-u_k(0)}\rightarrow 0,\quad\text{as}\,\, k\to+\infty,$$ such that, setting $$\eta_k(x):=u_k(\mu_k x)-u_k(0)+\log2,$$ we have $$\eta_k(x)\xrightarrow[k\rightarrow\infty]{}\log\frac{2}{1+|x|^2} \qquad\text{in}\quad C^1_{\text{loc}}(\mathbb{R}^2).$$
Since $u_k$ is monotone decreasing, we have that $u_k\rightarrow-\infty$ locally uniformly in $\mathbb{R}^2\setminus\{0\}$ and using \eqref{63} we get
\begin{equation}\label{65}
\lim\limits_{k\rightarrow+\infty}\int_{\mathbb{R}^2}e^{2u_k}dx=\Lambda_\text{sph}.\end{equation}
Comparing \eqref{65} and \eqref{62}, we have
$$1=\frac{2\bar{\Lambda}}{p\Lambda_\text{sph}}\ge\frac{2\Lambda_{*,p}}{p\Lambda_\text{sph}}=\frac{2+p}{2p}>1\qquad \text{for}\quad p\in(0,2),$$
which is absurd.
\endproof

\begin{lemma}\label{4.4}
$u_k\rightarrow\bar{u}$ locally uniformly, where $\bar{u}$ is a solution to \eqref{7} with $\Lambda=\tilde{\Lambda}\ge\bar{\Lambda}$.
\end{lemma}

\proof
Using Lemma \ref{4.2} and Lemma \ref{4.3} and the fact that $u_k$ is radial decreasing, we have $$u_k\le u_k(0)=O(1),$$
therefore $$-\Delta u_k=O_R(1) \qquad\text{on}\quad B_R.$$
Integrating \eqref{58} over $B_R$ and using Fubini's theorem and \eqref{63}
$$\int_{B_R}|u_k(x)|dx\le C_R$$
therefore, by elliptic estimates we have $$\|u_k\|_{\mathcal{C}^{1,\alpha}(B_{R/2})}\le C_{R/2}$$
and hence, up to a subsequence, $u_k\to\bar{u}$ in $\mathcal{C}^1_ {loc}(\mathbb{R}^2)$. Finally by Fatou's lemma we have
$$\tilde{\Lambda}:=\int_{\mathbb{R}^2}(1-|x|^p)e^{2\bar{u}} dx\ge \limsup\limits_{k\to+\infty}\int_{\mathbb{R}^2}(1-|x|^p)e^{2u_k} dx=\bar{\Lambda}.$$
\endproof

\begin{lemma}\label{4.5}
We have $\tilde{\Lambda}=\bar{\Lambda}$.
\end{lemma}

\proof
Assume by contradiction that $\tilde{\Lambda}>\bar{\Lambda}$. Using the fact that $u_k\to\bar{u}$ in $\mathcal{C}^0_{loc}(\mathbb{R}^2)$, we have
$$\tilde{\Lambda}-\bar{\Lambda}=-\int_{B_R^c}(1-|x|^p)e^{2u_k}dx+o_k(1)+o_R(1),\quad\text{as}\,\, k\to+\infty,\,\, R\to+\infty$$
hence
\begin{equation}\label{66}
0<\tilde{\Lambda}-\bar{\Lambda}=\lim\limits_{R\to+\infty}\lim\limits_{k\to+\infty}\int\limits_{B_R^c}(|x|^p-1)e^{2u_k}dx=\lim\limits_{R\to+\infty}\lim\limits_{k\to+\infty}\int\limits_{B_R^c}|x|^p e^{2u_k}dx=:\rho.
\end{equation}
We consider the Kelvin transform of $u_k$ as defined in \eqref{20}
$$\tilde{u}_k(x)=u_k\left(\frac{x}{|x|^2}\right)-\frac{\Lambda_k}{2\pi}\log|x|,\quad x\not=0.$$
From Proposition \ref{2.2} we have
$$\tilde{u}_k(x)=\frac{1}{2\pi}\int_{\mathbb{R}^2}\log\left(\frac{1}{|x-y|}\right)\,\left(1-\frac{1}{|y|^p}\right)\frac{e^{2\tilde{u}_k(y)}}{|y|^{2-p-\delta_k}}dy+c_k,\quad \delta_k=\frac{\Lambda_k-\Lambda_{*,p}}{\pi},$$
and with the same procedure for proving \eqref{30} we obtain
\begin{equation}\label{67}
\tilde{u}_k(x)=-\frac{1}{2\pi}\int\limits_{B_1}\log\left(\frac{1}{|x-y|}\right)\frac{e^{2\tilde{u}_k(y)}}{|y|^{2-\delta_k}}dy+O(1),\quad \text{if}\,\,\,0<|x|<1.
\end{equation}
If $\delta_k\not\to 0$, from \eqref{67} we have that $\tilde{u}_k=O(1)$ in $B_1$, but this contradicts the fact that $\rho>0$, hence we must have $\delta_k\to0$, namely $\Lambda_k\to\Lambda_{*,p}$.
Let us define $\varepsilon_k>0$ in such a way that
\begin{equation}\label{67a}
\int_{B_{\varepsilon_k}}\frac{e^{2\tilde{u}_k(y)}}{|y|^{2-\delta_k}}dy=\frac{\rho}{2},
\end{equation}
and we have that $\varepsilon_k\to0$ as $k\to\infty$. We can observe that for $y\in B_{\varepsilon_k}$ and $x\in B_{2\varepsilon_k}^c$\,\, $\log\left(\frac{1} {|x-y|} \right)=\log\left(\frac{1}{|x|}\right)+O(1)$, hence for $2\varepsilon_k\le|x|\le1$ we get
\begin{equation}\label{68}
\begin{split}
\tilde{u}_k(x)&=-\frac{\log(1/|x|)}{2\pi}\int\limits_{B_{\varepsilon_k}}\frac{e^{2\tilde{u}_k(y)}}{|y|^{2-\delta_k}}dy-\frac{1}{2\pi}\int\limits_{B_1\setminus B_{\varepsilon_k}}\log\left(\frac{1}{|x-y|}\right)\frac{e^{2\tilde{u}_k(y)}}{|y|^{2-\delta_k}}dy+O(1)\\
&\le-\frac{\rho}{4\pi}\log\left(\frac{1}{|x|}\right)+C.
\end{split}
\end{equation}
where in the last inequality, we used the fact that $\log\left(\frac{1}{|x-y|} \right)$ is lower bounded for $y\in B_1$ and $x\to0$.
From \eqref{68} we have 
\begin{equation}\label{69}
\lim\limits_{r\to0}\,\lim\limits_{k\to\infty}\sup_{B_r}\tilde{u}_k=-\infty
\end{equation}
and moreover 
$$\lim\limits_{r\to0}\,\lim\limits_{k\to\infty}\int_{B_r\setminus B_{2\varepsilon_k}}\frac{e^{2\tilde{u}_k(y)}}{|y|^{2-\delta_k}}dy=0,$$
which, using \eqref{66}, implies
\begin{equation}\label{69a}
\lim\limits_{k\to\infty}\int_{B_{2\varepsilon_k}}\frac{e^{2\tilde{u}_k(y)}}{|y|^{2-\delta_k}}dy=\rho.
\end{equation}
Hence from \eqref{67a}, \eqref{69} and \eqref{69a} we obtain
$$\frac{\rho}{2}=\lim\limits_{k\to\infty}\int_{B_{2\varepsilon_k}\setminus B_{\varepsilon_k}}\frac{e^{2\tilde{u}_k(y)}}{|y|^{2-\delta_k}}dy=o(1)\int_{B_{2\varepsilon_k}\setminus B_{\varepsilon_k}}\frac{1}{|y|^2}dy=o(1),\quad\text{as}\,\, k\to\infty,$$
which is absurd, because we had assumed $\rho>0$.
\endproof

\proof[Proof of Theorem \ref{teo1.3}]
Lemma \ref{4.4} and Lemma \ref{4.5} prove that, up to a subsequence, $u_k$ converges to $\bar{u}$ locally uniformly in $\mathbb{R}^2$, where $\bar{u}$ is a solution to \eqref{7} with $\Lambda=\bar{\Lambda}$. Estimate \eqref{10} follows from Lemma \ref{2.9}, while for \eqref{10a} we proceed in the same way as for proving \eqref{9}.
\endproof

\section{Proof of Theorem \ref{teo1.4}}

\begin{lemma}\label{5.1}
Let $(u_k)$ be a sequence of solutions to \eqref{7} with $\Lambda_k\uparrow \Lambda_\mathrm{sph}$. Then $$u_k(0)\to+\infty,\quad\text{as}\,\,k\to+\infty.$$
\end{lemma}

\proof
By Lemma \ref{4.2} we have that $u_k(0)\ge C$. Assume by contradiction that, up to a subsequence, $u_k(0)\to s\in \mathbb{R}$ as $k\to+\infty$. Proceeding as in the proof of Lemma \ref{4.4}, we get that $u_k\to\bar{u}$, where $\bar{u}$ is a solution to \eqref{7} for $\Lambda\ge4\pi$, which contradicts Theorem \ref{teo1.1}
\endproof

\proof[Proof of Theorem \ref{teo1.4}]
By Lemma \ref{5.1} the sequence $(u_k)$ blows up at the origin. Moreover, we have
$$\int_{B_1(0)}|\nabla u_k|dx\le C \int_{\mathbb{R}^2}\left(\int_{B_1(0)}\frac{1} {|x-y|}dx\right)(1+|y|^p)e^{2u_k(y)}dy\le C.$$
We can conclude using Theorem 2 in \cite{M}.
\endproof


\section{Appendix}

\proof[Proof of the estimate \eqref{stimau1}]
Since $u_1=u-u_2$ one has
\begin{align*}
u_1(x)&=\frac{1}{2\pi}\int\limits_{\mathbb{R}^2\setminus B_1(x)} \log\left(\frac{|y|}{|x-y|}\right)(1-|y|^p) e^{2u(y)}dy+\\ &\quad+\frac{1}{2\pi}\int\limits_{B_1(x)}\log\left(\frac{1}{|x-y|}\right)e^{2u(y)}dy+\frac{1}{2\pi}\int\limits_{B_1(x)}\log|y|(1-|y|^p)e^{2u(y)}dy+c=\\
&=\frac{1}{2\pi}\left(I_1+I_2+I_3\right)+c.
\end{align*}
First of all we estimate $I_1$ using the split of $\mathbb{R}^2$ introduced in the proof of Proposition \ref{2.3}
$$I_1=\int_{B_{R_0}(0)}(\cdot)dy+\int_{B_{\frac{|x|}{2}}\setminus B_1(x)}(\cdot)dy+\int_{\mathbb{R}^2\setminus (B_{|x|/2}\cup B_{R_0})}(\cdot)dy.$$
If $y\in B_{R_0}$ we have $\log\left(\frac{|y|}{|x-y|}\right)=-\log|x|+O(1)$ as $|x|\rightarrow\infty$, hence 
$$\int\limits_{B_{R_0}} \log\left(\frac{|y|}{|x-y|}\right)(1-|y|^p) e^{2u(y)}dy=(-\log|x|+O(1))(\Lambda+o_{R_0}(1)).$$
For $y\in B_{|x|/2}(x)\setminus B_1(x)$ we have $\log\left(\frac{|y|}{|x-y|}\right)=O(\log|x|)$ as $|x|\rightarrow\infty$, hence 
$$\int\limits_{B_{|x|/2}(x)\setminus B_1(x)} \log\left(\frac{|y|}{|x-y|}\right)(1-|y|^p) e^{2u(y)}dy=o_{R_0}(1)\,O(\log|x|).$$
For $y\in\mathbb{R}^2\setminus(B_{|x|/2}(x)\cup B_{R_0})$ we can observe that $\log\left(\frac{|y|}{|x-y|}\right)=O(1)$ as $|x|\to\infty$, so 
$$\int\limits_{\mathbb{R}^2\setminus(B_{|x|/2}(x)\cup B_{R_0})}\log\left(\frac{|y|} {|x-y|}\right)(1-|y|^p)e^{2u(y)}dy=O(1)\,o_{R_0}(1).$$
Concerning $I_2$, using \eqref{23} as $|x|\to\infty$ we have
$$I_2=\int_{B_1(x)}\log\left(\frac{1}{|x-y|}\right)e^{2u(y)}dy\le C\int_{B_1(x)}\log \left(\frac{1}{|x-y|}\right)|y|^{-\frac{\Lambda}{\pi}}dy$$
from Proposition \ref{lambda>0} $\Lambda\ge0$, hence if $\Lambda=0$ we have $I_2\le O(1)$ and if $\Lambda>0$ we obtain $I_2=o(1)$ as $|x|\to\infty$.
Let's consider $I_3$, as $|x|\to\infty$ we get
$$I_3=O\left(\log|x|\int_{B_1(x)}(1-|y|^p)e^{2u(y)}dy\right)=o(\log|x|)$$
where, since $(1-|y|^p)e^{2u(y)}\chi_{B_1(x)}\to 0$ a.e. as $|x|\to\infty$, using dominated convergence we have $\int_{B_1(x)}(1-|y|^p)e^{2u(y)}dy= \int_{\mathbb{R}^2}(1-|y|^p)e^{2u(y)}\chi_{B_1(x)}dy\to 0$ as $|x|\to\infty$.
Summing up, we can conclude that 
$$u_1(x)=-\frac{\Lambda+o(1)}{2\pi}\log|x|,\quad\mathrm{as}\,\,|x|\rightarrow+\infty.$$ 
\endproof

\end{document}